\documentclass[11pt]{article}
\usepackage{amssymb,amsfonts}

\newtheorem{theorem}{Theorem}[section]
\newtheorem{lemma}[theorem]{Lemma}
\newtheorem{proposition}[theorem]{Proposition}

\newtheorem{definition}[theorem]{Definition}
\newtheorem{corollary}[theorem]{Corollary}

\newenvironment{proof}{\bf Proof. \rm}{$\Box$}

\newcommand{\be}{\begin{equation}}
\newcommand{\ee}{\end{equation}}

\newcommand{\cA}{\mathcal{A}}
\newcommand{\cB}{\mathcal{B}}
\newcommand{\cC}{\mathcal{C}}
\newcommand{\cD}{\mathcal{D}}
\newcommand{\cE}{\mathcal{E}}
\newcommand{\cG}{\mathcal{G}}
\newcommand{\cL}{\mathcal{L}}
\newcommand{\cM}{\mathcal{M}}
\newcommand{\cF}{\mathcal{F}}
\newcommand{\cT}{\mathcal{T}}

\newcommand{\aA}[1]{ _{\alpha^{ #1 }}A}
\newcommand{\f}[1]{\phi_{\infty}^{( #1 )}}
\newcommand{\BDalf}{\cB_{\alpha}(\{n_k\})}
\newcommand{\BDthet}{\cB_{\theta}(\{n_k\})}
\newcommand{\dL}{\delta_{\mathbf{0}}}
\newcommand{\del}[1]{\delta_{\mathbf{#1}}}
\newcommand{\qnk}{Q(\delta(n_k))}

\begin{document}
\title{Limit C*-algebras associated with an automorphism}
\author{Baruch Solel \thanks{Supported by the Fund for the
Promotion of Research at the Technion}\\ Department of Mathematics
\\ Technion \\ Haifa 32000 \\ ISRAEL \\
\texttt{mabaruch@techunix.technion.ac.il}}

\maketitle

\begin{abstract}
We present and study C*-algebras generated by "periodic weighted
creation operators" on the Fock space associated with an
automorphism $\alpha$ on a C*-algebra $A$. These algebras can be
viewed as generalized Bunce-Deddens algebras associated with the
automorphism and can be written as a certain direct limit. We prove
 a crossed product presentation for such an
algebra and find a necessary and sufficient condition for it to be
simple. In the case where the automorphism is induced by an
irrational rotation (on $C(\mathbb{T})$) we compute the K-theory groups
 and obtain a complete
classification of these algebras.

\textbf{2000 Subject Classification} Primary: 46L05.
Secondary: 46L40, 46L35, 46L80

\textbf{Key words :} C*-algebras, Limit algebras, Bunce-Deddens
algebras, automorphism.

\end{abstract}

\begin{section}{Introduction}
The purpose of this paper is to to present and study a class of
C*-algebras generalizing the Bunce-Deddens algebras.

Recall that the Bunce-Deddens algebra $\cB(\{n_k\})$, associated
with the sequence $\{n_k\}$, is a quotient (by the compact
operators) of the C*-algebra generated by all weighted shifts
(with respect to a fixed basis) of period $n_k$ for some $k\geq
1$ (\cite{BD}).
 In fact, for a fixed $k\geq 1$, the C*-algebra of all weighted
shifts of period $n_k$ is isomorphic to the algebra of all
$n_k\times n_k$ matrices over the Toeplitz algebra and its
quotient by the compacts is isomorphic to
$M_{n_k}(C(\mathbb{T}))$. Thus the Bunce-Deddens algebra is a
direct limit of the algebras $M_{n_k}(C(\mathbb{T}))$.

A natural noncommutative generalization of the unilateral shift is
obtained by considering creation operators on Fock Hilbert space.
These operators were studied in applications to quantum physics and in free
probability. More recently, creation operators, and the algebras
that they generate, were studied in the context of Fock spaces of
C*-correspondences ( also called Hilbert C*-bimodules). This was
initiated by M. Pimsner (\cite{Pi}) and followed by several other
authors (e.g. \cite{MS},\cite{MS1},\cite{KPW},\cite{S} ).

In a recent paper, D. Kribs introduced the concept of periodic
weighted shifts on the full Fock space associated with a Hilbert space
of dimension $N$ (\cite{Kr}). The C*-algebra that he
obtained (taking the appropriate quotient of the C*-algebra generated
 by all periodic weighted
 shifts of period $n_k$ ) is a direct limit of algebras of the
 form $M_{d_k}(O_{N^k})$ where $d_k$ are certain integers and
 $O_{N^k}$ is the Cuntz algebra ( with $N^k$ generators).

In the present paper we start with an automorphism $\alpha$ on a
C*-algebra $A$ and form the (full) Fock space $\cF$ associated
with the pair $(A,\alpha)$ (i.e. associated with the bimodule
$_{\alpha}A$). As Pimsner showed in \cite{Pi},
 the C*-algebra
generated by the creation operators on this Fock space is the
Toeplitz extension ,$\cT(A,\alpha)$, of the crossed product
$A\times_{\alpha}\mathbb{Z}$. Letting $K(\cF)$ denote the compact
operators on $\cF$ (in the sense of compact operators on a
C*-module), the crossed product is the quotient of $\cT(A,\alpha)$
by $K(\cF)$. We shall define the concept of "periodic weighted
creation operators" (more precisely, we define , in Definition~\ref{weighted},
 weighted
representations of $_{\alpha}A$ on $\cF$) and consider the
C*-algebra generated by all these operators with a fixed period
$n$. This algebra turns out to be isomorphic to $M_n(\cT(A,\alpha^n))$
and, taking an appropriate quotient, we get
$M_n(A\times_{\alpha^n}\mathbb{Z})$.
Now let $\{n_k\}$ be an increasing sequence of positive integers
with $n_k|n_{k+1}$. Considering the C*-algebra generated by all
"periodic weighted creation operators" of period $n_k$ for some
$k\geq 1$ and taking an appropriate quotient we get a certain
direct limit of the algebras
$M_{n_k}(A\times_{\alpha^{n_k}}\mathbb{Z})$. We shall write
$\BDalf$ for this C*-algebra. It can be thought of as a
"generalized Bunce-Deddens algebra associated with $(A,\alpha)$".
The details of this construction are presented in
Section~\ref{prel}. We also show there that the algebra depends,
up to an isomorphism, on the supernatural number of $\{n_k\}$ (and
not on the sequence itself).

In Section~\ref{structure} we study the structure of the algebra
$\BDalf$. We present necessary and sufficient conditions for it
 to have a unique tracial state
(Theorem~\ref{trace}). In Theorem~\ref{isom} we show that the
algebra $\BDalf$ can also be written as a crossed product
$C(X,A)\times_{\sigma}\mathbb{Z}$. We then use it in
Theorem~\ref{simple} to prove a necessary and sufficient
condition for simplicity of the algebra.

In Section~\ref{th} we specialize to the case where
$A=C(\mathbb{T})$ and $\alpha$ is an irrational rotation by
$\theta$ ( so that $A\times_{\alpha}\mathbb{Z}$ is the irrational
rotation algebra $A_{\theta}$). In this case we use the notation $\BDthet$
for the algebra in order to
emphasize the dependence on $\theta$. It follows from the previous
section that the algebra is simple and has a unique tracial state
$\tau$. We show that the $K_0$ and $K_1$ groups of the algebra are
both
isomorphic to $Q(\delta)\oplus \mathbb{Z}$ (where $\delta$ is the
supernatural number of the sequence $\{n_k\}$ and $Q(\delta)$ is
the group of rational numbers that can be written as a quotient
$m/n_k$ for some $k$). We also find that the image of the map
$\tau_*$, defined on the $K_0$ group, is $Q(\delta)+\theta
\mathbb{Z}$. It then follows that two such algebras, with $\theta_1$,
$\theta_2$ and supernatural numbers $\delta_1$, $\delta_2$ are
isomorphic if and only if $\delta_1=\delta_2$ and either
$\theta_1+\theta_2$ or $\theta_1-\theta_2$ belongs to
$Q(\delta_1)$ (Corollary~\ref{isomorphism}).

\end{section}
\begin{section}{Preliminaries}\label{prel}
Let $A$ be a (unital) separable C*-algebra and let $\alpha$ be a
unital
automorphism of $A$.
Write $E= _{\alpha}A$ for the C*-correspondence associated with
this automorphism. As a vector space, $E$ is simply $A$. The bimodule structure is defined by
$$ a\cdot b \cdot c=\alpha(a)bc .$$
Here $a,c$ are in $A$ and $b$ is in $A$ viewed as an
element of $E$. The inner product is
$$\langle b_1,b_2 \rangle = b_1^*b_2 .$$
This makes $E$ into a C*-correspondence in the sense of \cite{MS}.
It is easy to check that $E^{\otimes k}= _{\alpha^k}A$ and the
isomorphism of $E^{\otimes k}\otimes E^{\otimes m}$ onto
$E^{\otimes (k+m)}$ is given by
$$ V_{k,m}(a\otimes b)= \alpha^m(a)b $$
for $a\in _{\alpha^k}A$ and $b\in _{\alpha^m}A$.
In what follows we tend to suppress $V_{k,m}$ and identify the two
spaces.

The Fock space associated with $_{\alpha}A = E$ is
$$\cF(E)=A \oplus E \oplus E^{\otimes 2} \oplus \ldots  .$$
As a C*-module it is just the direct sum of infinitely many copies
of $A$.
To distinguish between elements of $_{\alpha^k}A$ for different
$k$'s, we write $\delta_k$ for the unit of $A$ viewed as an
element of $_{\alpha^k}A$. Thus $\delta_k a$ would be $a$ as an
element of $\aA{k}$.
The Fock space has a left action that makes it into a
C*-correspondence:
$$\phi_{\infty}(a)\delta_kb=\delta_k\alpha^k(a)b. $$
We shall denote the Fock space associated with $\alpha^k$ by
$\cF^{(k)}$ and the left action there by $\f{k}$.

Every $a$ in $E$ defines a shift operator on $\cF$ ($=\cF (E)$) by
$$ T(a)\delta_k b=\delta_1a \otimes \delta_kb=
\delta_{k+1}\alpha^k(a)b.$$
(In the last equality we omitted the reference to $V_{1,k}$.)

A similar operator on $\cF^{(k)}$ will be denoted $T^{(k)}(a)$.

\begin{definition}\label{weighted}
Let $\lambda =\{\lambda_i : 1\leq i <\infty \}$ be a bounded sequence of
 positive elements
of $A$ . The \emph{weighted representation} of $E$ on
$\cF$ is
$$ T_{\lambda}(a)\delta_k b= \delta_{k+1}
\alpha^k(\lambda_{k+1}a)b $$
for $a,b \in A$.
The weighted representation will be said to be \emph{ periodic of period
$k$} if for all $i \geq 1$, $\lambda_{i+k}=\lambda_i $.
\end{definition}

\textbf{Note:}
the word "representation" above refers only to the fact that it is
a representation of $E$ as a right module over $A$ (i.e.
 $T_{\lambda}(a\cdot c)=T_{\lambda}(a)\phi_{\infty}(c)$ for $c\in A$ and
 $a\in E$).

\vspace{3mm}

Note that every operator on $\cF$ can be written as an infinite
matrix (with respect to the decomposition of $\cF$ as an infinite direct
 sum). It will be convenient to denote by $a\delta_{ij}$ the
operator that maps $\delta_jb$ to $\delta_iab$ (for $a$ in $A$).

Also we shall write $S$ for $T(1)$ (hence
$S\delta_ka=\delta_{k+1}a$).

\begin{definition}\label{Toeplitz}
The C*-algebra generated by the operators $\{\;\f{k}(A)\; ,\\
T^{(k)}(\aA{k}) \}$ on $\cF^{(k)}$ is called the Toeplitz algebra
 associated with
$\alpha^k$ and will be denoted $\cT_k$.
Since $\cF^{(k)}$ is a subspace of $\cF$, each $\cT_k$ can be
viewed as a subspace of $\cL(\cF)$.
\end{definition}

\vspace{3mm}

\begin{lemma}\label{algebra}
Given $k\geq 1$, the C*-algebra generated by all the periodic
representations of period $k$ (on $\cF$) is isomorphic to the
algebra of all $k\times k$ matrices over $\cT_k$.
\end{lemma}
\begin{proof}
We start by setting some notation. We can write
$$\cF=\cF_k \oplus (E\otimes \cF_k) \oplus (E^{\otimes 2}\otimes
\cF_k) \oplus \cdots \oplus (E^{\otimes (k-1)}\otimes \cF_k)$$
Let us write $\cM_l$ for $E^{\otimes l}\otimes \cF_k$ , $0\leq l
\leq k-1$. Then $\cF=\cM_0 \oplus \cM_1 \oplus \cdots \oplus
\cM_{k-1} $.
Every operator on $\cF$ can be written in a matricial way with
respect to this decomposition. We shall write $Xe_{ij}$ for the
operator whose $i,j$ entry is $X$ ($X$ maps $\cM_j$ into $\cM_i$)
and all other entries vanish.

Note that , for all $0\leq l \leq k-1$, $\phi_{\infty}(A)\cM_l
\subseteq \cM_l$ so that the matrix of $\phi_{\infty}(a)$ for
$a\in A$ is diagonal. For every $b\in E$, $T_{\lambda}(b)$ maps
$\cM_l$ into $\cM_{l+1}$ if $l<k-1$ and it maps $\cM_{k-1}$ into
$\cM_0$. Hence its matrix has non zero terms only on the first lower
diagonal and in the $0,k-1$ entry.

Now write
$$\cB=\sum_{i=0}^{k-1}\sum_{j=0}^{k-1} S^i\cT_k S^{*j}e_{ij} .$$
For $a\in A$ and $0\leq l \leq k-1$, the $l,l$ entry of
$\phi_{\infty}(a)$ is
$$\alpha^l(a)\delta_{ll}+\alpha^{l+k}(a)\delta_{l+k,l+k} +\cdots =
S^l (\alpha^l(a)\delta_{00}+\alpha^{k+l}(a)\delta_{kk}+\cdots )S^{*l}
=$$
$$=S^l\f{k}(\alpha^l(a))S^{*l} \in S^l
\cT_k S^{*l}. $$
For $b\in E$, the $l+1,l$ entry of $T_{\lambda}(b)$, for $l<k-1$, is
$$\alpha^l(\lambda_{l+1}b)
\delta_{l+1,l}+\alpha^{l+k}(\lambda_{l+1+k}b)\delta_{l+1+k,l+k}+\cdots. $$
 Since
$\lambda$ is of period $k$, this is equal to
$$\alpha^l(\lambda_{l+1}b)\delta_{l+1,l}+\alpha^k(\alpha^l(\lambda_{l+1}b))
\delta_{l+1+k,l+k}+\cdots $$
\be\label{l}
=S^{l+1}(\alpha^l(\lambda_{l+1}b)\delta_{00}+\alpha^{k+l}(\lambda_{l+1}b)
\delta_{k,k}+\cdots )S^{*l} \ee
$$= S^{l+1} \phi_{\infty}^{(k)}(\alpha^l(\lambda_{l+1}b))
S^{*l} \in S^{l+1}\cT_k S^{*l} .$$
For the $0,k-1$ entry of
$T_{\lambda}(b)$ we get
$$(
\alpha^{k-1}(\lambda_kb)\delta_{k,0}+\alpha^{2k-1}(\lambda_kb)\delta_{2k,k}+
\cdots )S^{*(k-1)} =$$
\be\label{k} =T^{(k)}(\alpha^{k-1}(\lambda_{k}b))S^{*(k-1)} \in
 \cT_k S^{*(k-1)} .\ee
Now write $\cA_k$ for the C*-algebra generated by all the periodic
representations of period $k$. It follows from the above that
$\cA_k \subseteq \cB $. We shall now show that, in fact, the two
algebras are equal.

Fix $0\leq j <k-1$. Setting $\lambda_i=1$ if $i
\equiv j+1 (mod\;k)$ and $\lambda_i=0$ otherwise, we find, using
 (\ref{l}) and (\ref{k}), that
$S^{j+1}\f{k}(\alpha^j(b))S^{*j}e_{j+1,j}=T_{\lambda}(b)\in \cA_k$. Since
this holds for all $b$ in $A$, we conclude that, for all $a\in
A$ and every $0\leq j <k-1$,
\be\label{a}
S^{j+1}\f{k}(a)S^{*j}e_{j+1,j} \in \cA_k .
\ee
Similarly, by setting $\lambda_i=1$ if $i\equiv 0\;(mod\;k)$ and
$\lambda_i=0$ otherwise
, we find that \be
T^{(k)}(E)S^{*(k-1)}e_{0,k-1} \subseteq \cA_k .
\ee
Thus, for every $b,a_1,\ldots a_{k-1}$ in $A$ we find that
$T^{(k)}(b)S^{*(k-1)}e_{0,k-1}$,\\ $
\;S^{k-1}\f{k}(a_{k-1})S^{*(k-2)}e_{k-1,k-2}$
, $\; \ldots$ and $ \;S\f{k}(a_1)e_{1,0} $ all lie in $\cA_k$. Multiplying
them, we see that $T^{(k)}(b)\f{k}(a_{k-1} \cdots a_1)e_{0,0}$ lie
in $\cA_k$.

 Thus $
T^{(k)}(E)e_{0,0} \subseteq \cA_k . $
Since $\{X : Xe_{0,0} \in \cA_k \} $ is a C*-algebra, \be\label{00}
\cT_k e_{0,0} \subseteq \cA_k . \ee
Setting, in (\ref{a}), $a=1$ we find that
$S^{j+1}S^{*j}e_{j+1,j}$ lies in $\cA_k$ (for $0\leq j <k-1$).
 It follows that
  $$S^{j+1}e_{j+1,0}=(S^{j+1}S^{*j}e_{j+1,j})(S^jS^{*(j-1)}e_{j,j-1}) \cdots
   Se_{1,0}\in \cA_k.$$ Finally, for
 $0\leq i,j \leq k-1$ and $X$ in $\cT_k$,
$$ S^jXS^{*i}e_{i,j}= (S^je_{j,0})(Xe_{0,0})(S^ie_{i,0})^* \in
\cA_k. $$
Hence $\cB = \cA_k$.
But $\cB$ is clearly isomorphic to $M_k(\cT_k)$ and the
isomorphism is given by conjugating by the (isometric) matrix
 $diag(1,S,S^2,\ldots, S^{k-1})$.
\end{proof}

\vspace{3mm}

Suppose $n,m$ are two natural numbers with $n|m$ and let $\lambda$
be periodic of period $n$. Then $\lambda$ is periodic of period
$m$ and this implies that $\cA_n \subseteq \cA_m$. Using the
isomorphisms of Lemma~\ref{algebra} we get an injective
*-homomorphism $\beta_{n,m} :M_n(\cT_n) \rightarrow M_m(\cT_m) $. In the
next result we describe this map. In what follows we keep the
notation set up in the proof of Lemma~\ref{algebra}.

\vspace{3mm}

\begin{lemma}
Let $n,m$ be two natural numbers with $n|m$. Write $k=m/n$. Then
the restriction of $\beta_{m,n}$ to $\cT_n e_{0,0}$ is defined by
$$ \beta_{n,m}(\f{n}(a)e_{0,0})= \sum_{j=0}^{k-1}
\f{m}(\alpha^{jn}(a))e_{jn,jn} $$
for $a\in A$ and
$$\beta_{n,m}(T^{(n)}(b)e_{0,0})=\sum_{j=0}^{k-2}
\f{m}(\alpha^{jn}(b))e_{(j+1)n,jn} +
T^{(m)}(\alpha^{(k-1)n}(b))e_{0,(k-1)n} $$
for $b\in E$.

Also, for $0\leq i,j \leq n-1$ we have
\be\label{ij} \beta_{n,m}(e_{i,j})=\sum_{l=0}^{k-1} e_{i+ln,j+ln} .\ee
\end{lemma}
\begin{proof}
In the following we shall write $\eta_n$ for the isomorphism
$\eta_n: \cA_n\rightarrow M_n(\cT_n)$ of Lemma~\ref{algebra}.
Fix $0\leq j <n-1$ and let $\lambda_p$ be $1$ if $p\equiv j+1 (mod
\;n)$ and be $0$ otherwise. Then (by equations (\ref{l})
and (\ref{k})) the matricial form of $T_{\lambda}(1)$ as an element
of $\cA_n$ is $S^{j+1}S^{*j}e_{j+1,j}$. Thus
$\eta_n(T_{\lambda}(1))=e_{j+1,j}$.
 Using  (\ref{l})
and (\ref{k}) again, this time with $m$ instead of $n$, we find
that the matricial form of $T_{\lambda}(1)$ as an element of
$\cA_m$ is $\sum_{l=0}^{k-1} S^{j+nl+1}S^{*(j+nl)}e_{j+nl+1,j+nl}$
and, thus, $\eta_m(T_{\lambda}(1))=\sum e_{j+nl+1,j+nl}$. Equation
(\ref{ij}) for $i=j+1$ follows. Since $\beta_{n,m}$ is a
*-homomorphism, (\ref{ij}) follows for all $i,j$.

Now let $\lambda$ be defined by $\lambda_p=1$ if $p$ is a multiple
of $n$ and is equal to $0$ otherwise. Then, for $b\in E$, the matricial
 form of $T_{\lambda}(\alpha^{1-n}(b))$ as an element of $\cA_n$
 is $T^{(n)}(b)S^{*(n-1)}e_{0,n-1}$ (hence $\eta_n(T_{\lambda}(\alpha^{1-n}(b)))
 =T^{(n)}(b)e_{0,n-1}$ ) and its matricial form as an
 element of $\cA_m$ is
 $T^{(m)}(\alpha^{(k-1)n}(b))S^{*(m-1)}e_{0,m-1}
 +\sum_{l=1}^{k-1}S^{ln}\f{m}(\alpha^{l-1}(b))S^{*(ln-1)}e_{ln,ln-1}$
 (and therefore $\eta_m(T_{\lambda}(\alpha^{1-n}(b)))=\\ T^{(m)}(\alpha^{(k-1)n}(b))e_{0,m-1}+
\sum_{l=1}^{k-1} \f{m}(\alpha^{l-1}(b))e_{ln,ln-1}$).
 It follows that
$$\beta_{n,m}(T^{(n)}(b)e_{0,n-1})=T^{(m)}(\alpha^{(k-1)n}(b))e_{0,m-1}+
\sum_{l=1}^{k-1} \f{m}(\alpha^{l-1}(b))e_{ln,ln-1} .$$
Since
$$T^{(n)}(b)e_{0,0}=(T^{(n)}(b)e_{0,n-1})\cdot e_{n-1,n-2} \cdots e_{1,0}
$$
and,
$$\f{n}(a)e_{0,0}=(T^{(n)}(1)e_{0,0})^*T^{(n)}(a)e_{0,0}, \;\;\;a\in A, $$
a straightforward computation, using the fact that $\beta_{n,m}$
is a *-homomorphism, completes the proof.

\end{proof}

\begin{corollary}\label{beta}
For $n|m$ the map $\beta_{n,m}$, described above, maps $\cT_n
e_{i,j}$ into $\sum_{l=0}^{k-1} \cT_m e_{i+ln,j+ln} $ where
$k=m/n$.
\end{corollary}

\vspace{3mm}

In fact, define the map $\theta_{n,m} :\cT_n \rightarrow
M_k(\cT_m) $ by
$$
\theta_{n,m}(\f{n}(a))=\sum_{j=0}^{k-1}\f{m}(\alpha^{jn}(a))e_{j,j}$$
for $a\in A$ and
$$ \theta_{n,m}(T^{(n)}(b))=T^{(m)}(\alpha^{(k-1)n}(b))e_{0,k-1} +
\sum_{j=0}^{k-2}\f{m}(\alpha^j(b))e_{j+1,j} $$
for $b\in E$.

The matricial form of the elements of $\cA_m$ was with respect to
the decomposition
$$\cF=\cM_0 \oplus \cdots \oplus \cM_{m-1} $$
where $\cM_l=E^l \otimes \cF_m$. As a C*-module over $A$ each
$\cM_l$ is isomorphic to $\cF_m$. Thus any permutation of these
spaces defines a unitary operator on $\cF$. Let $U$ be the
operator arising from the permutation $(0,1,\ldots ,m-1) \mapsto
(0,n,\ldots ,(k-1)n,1,2,\ldots ,m-n-1,m-1)$. Then we have the
following.

\begin{corollary}\label{shuffle}
With the notation above,
$$ \beta_{n,m}=U^*(I_k \otimes \theta_{n,m})U .$$
\end{corollary}

\vspace{3mm}

The following lemma follows immediately from the definition
 of the maps $\beta_{n,m}$ (or by a straightforward tedious
computation using Lemma~\ref{beta}).

\begin{lemma}\label{composition}
For positive integers $n,k,l$ we have
$$ \beta_{nk,nkl}\circ \beta_{n,nk}=\beta_{n,nkl} .$$
\end{lemma}

\vspace{3mm}

It is known that the Toeplitz algebra $\cT_n$ is an extension of
the crossed product algebra of $A$ by $\alpha^n$, $A
\times_{\alpha^n} \mathbb{Z} $. The Toeplitz algebra
contains the algebra $K(\cF_n)$ of all compact operators on the
C*-module $\cF_n$ (where "compact" is in the sense of C*-modules
theory ) which is isomorphic to $A\otimes K$ (where $K$ denotes
the algebra of compact operators on a separable Hilbert space).
The quotient space $\cT_n / K(\cF_n)$ is then isomorphic to the
crossed product. One can also check (see e.g. \cite{MS}) that the
ideal $K(\cF_n)$ is generated by the set $\f{n}(A)P_0$ where $P_0$
is the projection of $\cF_n$ onto the first summand of $\cF_n$ (i.e. $A$) .

We now compute, for $a,b,c$ in $A$ and $k\geq 1$,
$$(\phi_{\infty}(ac^*)-T(\alpha(a))T(\alpha(c))^*)\delta_kb
=\delta_k\alpha^k(ac^*)b-T(\alpha(a))\delta_{k-1}\alpha^k(c^*)b
=$$
$$=\delta_k\alpha^k(ac^*)b-\delta_k \alpha^k(ac^*)b=0 $$
and, for $k=0$,
$$(\phi_{\infty}(ac^*)-T(\alpha(a))T(\alpha(c))^*)\delta_0b =
\delta_0ac^*b .$$
Hence,
\be\label{id} \phi_{\infty}(ac^*)-T(\alpha(a))T(\alpha(c))^*=
\phi_{\infty}(ac^*)P_0.   \ee

\begin{lemma}\label{cp}
For every positive integers $n,m$ with $ m/n =k\in \mathbb{Z}$ the
map $\beta_{n,m}$ induces a map, denoted $\gamma_{n,m}$, from
$M_n(A\times_{\alpha^n} \mathbb{Z}) $ into $M_m(A\times_{\alpha^m}
\mathbb{Z})$.
\end{lemma}
\begin{proof}
Using Corollary~\ref{shuffle} it suffices to show that a similar
statement holds for the map $\theta_{n,m}$.
For this, we have to show that $\theta_{n,m} $ maps $K(\cF_n)$ into
$M_k(K(\cF_m))$.
We compute, for $a,c\in A$,
$$\theta_{n,m}(\f{n}(ac^*)-T^{(n)}(\alpha^n(a))T^{(n)}(\alpha^n(c))^*)=$$
$$=\sum_{j=0}^{k-1}
\f{m}(\alpha^{jn}(ac^*))e_{j,j}-T^{(m)}(\alpha^{(k-1)n}(\alpha^n(a)))T^{(m)}
(\alpha^{(k-1)n}(\alpha^n(c)))^*e_{0,0}$$
$$-\sum_{j=0}^{k-2}
 \f{m}(\alpha^{jn}(\alpha^n(ac^*)))e_{j+1,j+1}=$$
$$=( \f{m}(ac^*)-T^{(m)}(\alpha^m(a))T^{(m)}(\alpha^m(c))^*)e_{0,0}. $$
From the discussion preceeding the lemma it now follows that
$\theta_{n,m}$ maps the set $\f{n}(A)P_0$ into $M_k(K(\cF_m))$.
Since this set generates the ideal $K(\cF_n)$, the claim follows.

\end{proof}

\vspace{3mm}

Let $\gamma_{m,n}$ be the map defined above and let $q_n$ and $q_m$ be the
quotient maps on $M_n(\cT_n)$ and $M_m(\cT_m)$ respectively
(induced by $q$). Recall that $q(\f{n}(a))=a$ for $a\in A$ and
$q(T^{(n)}(1))=u_n^*$ ( where $u_n$ is the unitary that, together
with $A$, generates $A\times_{\alpha^n}\mathbb{Z}$ and satisfies $u_nau_n^*=
\alpha^n(a)$, $a\in A$). Then it follows
from the definition of $\beta_{n,m}$ that $\gamma_{n,m}$ is
defined by
$$\gamma_{n,m}(ae_{0,0})=\sum_{l=0}^{k-1} \alpha^{ln}(a)e_{ln,ln}
$$ for $a\in A$,
$$\gamma_{n,m}(u_n)=u_me_{(k-1)n,0}+\sum_{l=0}^{k-2}e_{ln,(l+1)n}
$$ and, for $0\leq i,j \leq n-1$,
$$\gamma_{n,m}(e_{i,j})=\sum_{l=0}^{k-1} e_{i+ln,j+ln}. $$

We now write $\cB(n)$ (or $\cB_{\alpha}(n)$ if we want to emphasize
 the dependence on $\alpha$) for the algebra $M_n(A\times_{\alpha^n}
\mathbb{Z})$.

\begin{definition}\label{Balpha}
 Suppose that $\{n_k\}$ is an increasing sequence of
positive integers such that $n_k$ divides $n_{k+1}$ for $k\geq 1$.
The \emph{Bunce-Deddens algebra of $\alpha$} is the direct limit
$$ \cB_{\alpha}(\{n_k\})=lim_{\rightarrow} (\cB(n_k), \gamma_{n_k,n_{k+1}})
.$$
\end{definition}

Recall that an element $\sum_{k=-\infty}^{\infty}a_ku_n^k$ in the
crossed product is $0$ if and only if $a_k=0$ for every $k$. It
follows from this and the definition of $\gamma_{n,m}$ that this map
is injective . We can, thus, write
the limit algebra $\cB_{\alpha}(\{n_k\})$ as the closure of the increasing
union $\cup \cB(n_k) $.

Note that the analysis above shows that, for every $k \geq 1$, we
have
$$\cB(n_k)=M_{n_k}(A\times_{\alpha^{n_k}}\mathbb{Z}) \cong
M_{n_k}(\cT_{n_k})/ M_{n_k}(K(\cF_{n_k})) \cong \cA_{n_k}/
K(\cF). $$
It follows that $\BDalf$ is isomorphic to the quotient of
$\overline{\cup\cA_{n_k}}$ by $K(\cF)$. Recall that $\overline{\cup
\cA_{n_k}}$ is the C*-algebras generated by all periodic
representations of perion $n_k$ for all $k\geq 1$.

\begin{proposition}\label{class}
Let $\{n_k\}$ and $\{m_j\}$ be increasing sequences of positive
integers for which $n_k|n_{k+1}$ and $m_j|m_{j+1}$ for all $j,k
\geq 1$. Write $\delta(n_k)$ and $\delta(m_j)$ for the
corresponding supernatural numbers. If
$\delta(n_k)=\delta(m_j)$ then $\cB_{\alpha}(\{n_k\})$ and
 $\cB_{\alpha}(\{m_j\})$
are isomorphic.
\end{proposition}
\begin{proof}
The proof is standard. For each $k$, $n_k$ divides some $m_r$ and
the map $\gamma_{n_k,m_r}$ is the inclusion of $\cB(n_k)$ into
$\cB(m_r)$. Since the latter algebra is contained in
$\cB(\{m_j\})$, we find that every $\cB(n_k)$ is contained there.
Since this holds for all $k$, $\cB(\{n_k\}) \subseteq
\cB(\{m_j\})$. Equality follows from symmetry.
\end{proof}

\vspace{3mm}

the converse of the proposition above holds in some cases (see
Corollary~\ref{isomorphism}) but not in general, as we see in the
following lemma.

\begin{lemma}\label{converse}
Suppose $A$ is a unital C*-algebra such that $A$ is isomorphic to
$M_p(A)$ for some positive integer $p$. Let $\{n_k\}$ be an
increasing sequence as above and let $\alpha=id$ (the identity
automorphism). Then
$$ \cB_{id}(\{n_k\}) \cong \cB_{id}(\{pn_k\}) . $$
\end{lemma}
\begin{proof}
Since $\alpha=id$, $A\times_{\alpha^n}\mathbb{Z} \cong A\otimes
C(\mathbb{Z})$ for all $n$. Making these identifications, the maps
$\gamma_{n_k,n_{k+1}}$ (for the algebra $\cB_{id}(\{n_k\})$
 now satisfy (when we write $m_k$ for
$n_{k+1}/n_k$)
$$\gamma_{n_k,n_{k+1}}((a\otimes 1)e_{i,j})=\sum_{l=0}^{m_k-1}
(a\otimes 1)e_{i+ln_k,j+ln_k} $$
for $a\in A$, and
$$\gamma_{n_k,n_{k+1}}((a\otimes z)e_{i,j})=(a\otimes
z)e_{i+(m_k-1)n_k,j} +\sum_{l=0}^{m_k-2}(a\otimes 1) e_{i+n_kl,j+n_k(l+1)} $$
where $z$ is the identity function on $\mathbb{T}$.
The maps $\gamma_{pn_k,pn_{k+1}}$ (for the other algebra) can be
written similarly.
Write $\psi$ for the isomorphism $\psi:A\rightarrow M_p(A)$. It
induces an isomorphism from $A\otimes C(\mathbb{T})$ onto
$M_p(A\otimes C(\mathbb{T}))$ defined by
$$t(a\otimes f)=\sum_{0\leq i,j \leq p-1} (\psi(a)_{i,j} \otimes
f)e_{i,j} .$$
For every positive integer $k$ we write $t_k$ for the isomorphism
from $M_k(A\otimes C(\mathbb{T}))$ onto $M_{kp}(A\otimes
C(\mathbb{T}))$ defined by applying $t$ to each entry. Hence, for
$0\leq r,q \leq k-1$,
$$t_k((a\otimes f)e_{r,q})=\sum_{0\leq i,j \leq p-1}
(\psi(a)_{i,j} \otimes f)e_{i+rp,j+qp} .$$
To prove the claimed isomorphism it suffices to show that, for
every $k \geq 1$,
\be\label{t}
\gamma_{pn_k,pn_{k+1}} \circ t_{n_k} = t_{n_{k+1}} \circ
\gamma_{n_k,n_{k+1}} .
\ee
The required isomorphism in this case is
$\lim_{\rightarrow}t_{n_k}$.
To prove (\ref{t}) it is enough to apply both sides to elements of
the form $(a\otimes z)e_{0,0}$ and $(a\otimes 1)e_{r,q}$ for $a\in
A$. The computation is straightforward and is omitted.
\end{proof}

\end{section}

\begin{section}{The structure of $\BDalf$}\label{structure}

\begin{proposition}\label{trace}
$\BDalf$ has a unique tracial state $\tau$ if and only if $A$ has an
$\alpha$-invariant tracial state $\tau_0$ and for every $k\geq 1$
this is the unique $\alpha^{n_k}$-invariant tracial state on $A$.

Moreover, in this case, $\tau$ is faithful if and only if $\tau_0$
is.

\end{proposition}

\begin{proof}
Let $\tau_0$ be an $\alpha^{n_l}$-invariant tracial state on $A$.
For every $k\geq 1$, let
$\cE_k $ be the canonical conditional expectation of the crossed
product $A\times_{\alpha^{n_k}}\mathbb{Z}$ onto $A$ and write
$\tau_k=\tau_0\circ \cE_k$. If $k\geq l$ then $\tau_k$ is a tracial state on the
crossed product ( faithful if $\tau_0$ is). For such $k$ and an element
$X=(x_{ij})$ of the $n_k\times n_k$ matrices over the crossed
product $A\times_{\alpha}^{n_k}\mathbb{Z}$, we set
$\tilde{\tau}_k(X)=1/n_k (\sum \tau_k(x_{ii}))$ and this
defines a tracial state on $\cB(n_k)$ (that is faithful if $\tau_0$
is). To show that these traces define a tracial state on the direct
 limit it is left to show that for
every $k\geq l$,
\be \label{tau}
\tilde{\tau}_{k+1}\circ
\beta_{n_k,n_{k+1}}=\tilde{\tau}_{k} .
\ee
For $X=ae_{0,0}$ we have (set $m=n_{k+1}/n_k$),
$$\tilde{\tau}_{k+1}\circ \gamma_{n_k,n_{k+1}}(ae_{0,0})=
\tilde{\tau}_{k+1}(\sum_{j=0}^{m-1} \alpha^{jn_k}(a)e_{jn_k,jn_k})
= \sum
1/n_{k+1}\tau_k(\alpha^{jn_k}(a))=$$
$$=1/n_k\tau_0(a)=\tilde{\tau}_k(ae_{0,0}) .$$
For $X=u_{n_k}e_{0,0}$ we have,
$$\tilde{\tau}_{k+1} \circ
\gamma_{n_k,n_{k+1}}(u_{n_k}e_{0,0})=\tilde{\tau}_{k+1}(u_{n_{k+1}}e_{0,(k-1)n_k}
+\sum \tilde{\tau}_{k+1}(e_{(j+1)n_k,jn_k}) =0=$$
$$=\tau_k(u_{n_k})=\tilde{\tau}_k(u_{n_k}e_{0,0}) .$$
Since every tracial state on
$M_{n_k}(A\times_{\alpha^{n_k}}\mathbb{Z})$ is determined by its
values on elements of the form $Xe_{0,0}$ ($X\in
A\times_{\alpha^{n_k}}\mathbb{Z}$), this proves (\ref{tau}).
Thus, we constructed a tracial state $\tau$ on $\BDalf$ that is
faithful if $\tau_0$ is. Write $\tau=\Phi_l(\tau_0)$. So $\Phi_l$
maps $\alpha^{n_l}$-invariant traces of $A$ to traces of $\BDalf$.
It follows from the definition of $\Phi_l$ that it is one-to-one.

To go in the other direction, fix a tracial state $\tau$ on
$\BDalf$ and define, for $y$ in $A\times_{\alpha^{n_k}}\mathbb{Z}$
, $\tau_k(y)=1/n_k \tau(ye_{0,0})$ (where $e_{0,0}$ is an
 element of $\cB(n_k)$). Note that we have the following.
 $$\tau_k(u_{n_k})=1/n_k\tau(u_{n_k}e_{0,0})=1/n_k\tau(\gamma_{n_k,n_{k+1}}
 (u_{n_k}e_{0,0}))=$$
 $$=1/n_k
 (\tau(u_{n_{k+1}}e_{n_{k+1}-n_k,0})-\sum_{l=0}^{m-2}\tau(e_{ln_k,(l+1)n_k}))
 =0 $$
 where $m=n_{k+1}/n_k$. The last term is equal to $0$ because we
 evaluate $\tau$ there on matrices whose diagonal entries are all
 zeros. It follows that $\tau_k=(\tau_k|A) \circ \cE_k$ ; i.e.
 $\tau_k$ is determined by its values on $A$. Now set $\Psi_k(\tau)=
\tau_k| A$ for $a\in A$. Since $\tau_k$ is a tracial state on
$A\times_{\alpha^{n_k}}\mathbb{Z}$, it follows that $\Psi_k(\tau)$ is
$\alpha^{n_k}$-invariant.

Now suppose that $\BDalf$ has a unique tracial state $\tau$. Let
$\tau_0$ be $\Psi_1(\tau)$. Then it is an $\alpha$-invariant tracial
state on $A$. Fix $k\geq 1$. Suppose $\tau'_0$ is another
(different) $\alpha^{n_k}$-invariant trace on $A$. Then
$\Phi_k(\tau'_0)\neq \Phi_k(\tau_0)$ (since $\Phi_k$ is injective).
But this contradicts the assumed uniqueness of traces on $\BDalf$.
This proves one direction.
Note also that if $\tau$ is faithful then so is $\tau_0$ (which
is, roughly speaking, a restriction of $\tau$ to a copy of $A$).

 For the other direction, let
$\tau_0$ be an $\alpha$-invariant trace with the uniquness
property stated in the proposition. Write $\tau=\Phi_1(\tau_0)$.
Then $\tau$ is a tracial state on $\BDalf$. Suppose $\tau'$ ia
also a tracial state on $\BDalf$ that is different from $\tau$.
But then, for some $k\geq 1$, the restrictions of $\tau$ and
$\tau'$ to $\cB(n_k)$ are different. Since a trace on
$\cB(n_k)=M_{n_k}(A\times_{\alpha^{n_k}}\mathbb{Z})$ is determined
by its values on matrices of the form $ye_{0,0}$ for $y\in
A\times_{\alpha^{n_k}}\mathbb{Z}$, and since, as we saw above,
these values are determined by the restriction to elements of the
form $ae_{0,0}$ for $a\in A$,
 we find that $\Psi_k(\tau)\neq
\Psi_k(\tau')$. But this contradicts the uniqueness property of
$\tau_0$. In the construction of $\Phi_1(\tau_0)$ above it was
noted that if $\tau_0$ is faithful so is $\Phi_1(\tau_0)$.

\end{proof}

\vspace{3mm}

We shall now proceed to show that the algebra $\BDalf$ is
isomorphic to a crossed product algebra of the form
$C(X,A)\times_{\sigma}\mathbb{Z}$.

In order to define the space $X$, fix a sequence $\{n_k\}$ of
positive integers with $n_k|n_{k+1}$ for all $k\geq 1$ and with
$n_1=1$. Set $m_k=n_{k+1}/n_k$, let $X_i=\{0,1, \ldots ,m_i-1\}$
 and let $X$ be the Cantor set
$X=\prod_{i=1}^{\infty}X_i$. We can think of each element $x$ of
$X$ either as a sequence $\{x_i\}$ or as a formal sum $\sum_i
x_in_i$. We write $\sigma_0$ for the odometer action on $X$. For
$x\in X$ let $i(x)$ be the smallest $i$ for which $x_i<m_i-1$
(and $i(x)=\infty$ if $x_i=m_i-1$ for all $i$). Then, if
$i(x)<\infty$, $\sigma_0(x)_i=0$ for all $i< i(x)$ ,
$\sigma_0(x)_{i(x)}=x_{i(x)}+1$ and $\sigma_0(x)_i=x_i$ for $i>i(x)$.
 If $i(x)=\infty$, $\sigma_0(x)=\{0,0,
\ldots \}$. It is known that this map is a homoemorphism on $X$.

We consider the cylinder sets
$$J(x_1,x_2,\ldots x_k)=\{y \in X : y_i=x_i \;\; 1\leq i \leq k
\}.$$
Then $\sigma_0(J(x_1,\ldots,x_k))=J(y_1,\ldots,y_k)$ where
$$1+\sum_{i=1}^kx_in_i \equiv \sum_{i=1}^ky_in_i \;\;\; (mod\;
n_{k+1}). $$
Given  $(y_1,\ldots,y_k)$ with $0\leq
y_i \leq m_i-1$, we write
$\delta_{(y_1,\ldots,y_k)}$ for the function in
$C(X,A)$ which is $1$ (the unit of $A$) on $J(y_1,\ldots,y_k) $
and $0$ otherwise.
For an element $a\in A$ we shall write (by a slight abuse of
notation) also $a$ for the function in $C(X,A)$ that is constantly
equal to $a$. For $k\geq 2$ write
$$ \cC_k=\{f \in C(X,A) : f(x)=f(y) \;\; whenever
\;\;(x_1,\ldots,x_{k-1})=(y_1,\ldots,y_{k-1}) \} $$
and, for $k=1$, $\cC_1=A$ (where $A$ is viewed as the algebra of
all constant functions in $C(X,A)$).
Then $\cC_k \subseteq \cC_{k+1}$ and
$$ C(X,A)=\overline{\cup_{k=1}^{\infty} \cC_k }.$$
Now define the automorphism $\sigma$ on $C(X,A)$ by
$$\sigma(f)(x)=\alpha(f(\sigma_0^{-1}(x))) .$$
We shall show that $\BDalf$ is isomorphic to
$C(X,A)\times_{\sigma}\mathbb{Z}$. To do this, we write $U$ for
the unitary that , together with $C(X,A)$, generates the crossed
 product and let, for
every $k\geq 1$, $\cG_k$ be the C*-algebra generated by $\cC_k$
and $U$. For a given $k$ and $0\leq j\leq n_k-1$, write $\mathbf{j}$ for
$(j_1,\ldots,j_{k-1})$ satisfying $\sum_l j_ln_l=j$. Note that we have
$$ U^j\delta_{\mathbf{0}}U^{*j}=\sigma^j(\delta_{\mathbf{0}})=
\delta_{\mathbf{0}}\circ \sigma_0^{-j}=\delta_{\mathbf{j}} $$
for $0\leq j \leq n_k-1$. For $j=n_k$ a similar argument shows
that
$$ U^{n_k}\delta_{\mathbf{0}}U^{*n_k}=\delta_{\mathbf{0}} .$$

\begin{lemma}\label{Gk}
The algebra $\cG_k$ is the C*-algebra generated by the algebra $A$
(as a subalgebra of $\cC_k$), the operator
$\delta_{\mathbf{0}}U^{n_k}\delta_{\mathbf{0}}$ and the operators
of the form $U^i\delta_{\mathbf{0}}U^{*j}$ for $0\leq i,j \leq
n_k-1$. (For $k=1$ we let $\dL$ be $1$).
\end{lemma}
\begin{proof}
Write $\cD_k$ for the algebra defined in the statement of the
lemma. Since the function $\delta_{\mathbf{0}}$ lies in $\cC_k$,
it follows that $\cD_k\subseteq \cG_k$. Given $f\in \cC_k$ and
$0\leq j \leq n_k-1$ we write $a_j=f(j_1,j_2,\ldots )$ where
$j=\sum_{l=1}^{k-1} j_ln_l$ . Then
$$f=\sum_{j=0}^{n_k-1} a_j\delta_{\mathbf{j}}=\sum_j
a_jU^j\delta_{\mathbf{0}}U^{*j}.$$
Thus $\cC_k\subseteq \cD_k$. Also
$$U=\sum_{j=0}^{n_k-1}U\delta_{\mathbf{j}}=\sum_{j=0}^{n_k-2}U(U^j
\delta_{\mathbf{0}}U^{j*})+U(U^{n_k-1}\delta_{\mathbf{0}}U^{(n_k-1)*})=$$
$$=\sum_{j=0}^{n_k-2}U^{j+1}\delta_{\mathbf{0}}U^{j*}+
(\delta_{\mathbf{0}}U^{n_k}\delta_{\mathbf{0}})
(\delta_{\mathbf{0}}U^{(n_k-1)*}).$$
Thus $U\in \cD_k$, completing the proof.
\end{proof}

\begin{lemma}\label{rho}
Fix $k\geq 1$ and define the map
$$\rho_k:M_{n_k}(A\times_{\alpha^{n_k}}\mathbb{Z}) \rightarrow
\cG_k$$
by
$$
\rho_k(au_{n_k}^le_{i,j})=U^{*i}a\dL U^{j+n_kl}.$$
Then $\rho_k$ is a *-isomorphism onto $\cG_k$.
\end{lemma}
\begin{proof}
Note first that, since $a\delta_{\mathbf{0}}$ belongs to $\cC_k$,
the image of $\rho_k$ lies in $\cG_k$.
Write $x=au_{n_k}^le_{i,j}$ and $y=bu_{n_k}^re_{p,q}$ for $a,b\in
A$, $r,l \in \mathbb{Z}$ and $0\leq i,j,q,p \leq n_k-1$.
We compute (assuming for simplicity that $j \geq p$)
$$\rho_k(x)\rho_k(y)=U^{*i}a\dL U^{n_kl+j}U^{*p}b\dL
U^{n_kr+q}=$$
$$=U^{*i}a\dL \alpha^{j-p+n_kl}(b)\del{j-p}U^{n_k(l+r)-p+j+q}
.$$
 If $j\neq
p$, $xy=0$ and the computation above shows that also $\rho_k(x)\rho_k(y)
=0$ because $\dL \del{j-p}=0$ in this case. Suppose $p=j$. Then the
 computation above shows that
$$\rho_k(x)\rho_k(y)=U^{*i}a\alpha^{n_kl}(b)\dL
U^{n_k(l+r)}U^{q}=\rho_k(a\alpha^{n_kl}(b)u_{n_k}^{l+r}e_{i,q})=$$
$$=\rho_k(au_{n_k}^lbu_{n_k}^re_{i,q})=\rho_k(xy) .$$
Thus $\rho_k$ is multiplicative.
To show that $\rho_k$ is a *-map we compute
$$(\rho_k(au_{n_k}^le_{i,j}))^*=U^{-j-n_kl}\dL a^*
U^i=U^{*j}\alpha^{-n_kl}(a^*)\dL
U^{i-n_kl}=$$
$$=\rho_k(\alpha^{n_kl}(a^*)u_{n_k}^{*l}e_{j,i})=
\rho_k((au_{n_k}^le_{i,j})^*) .$$

Note that
$$\del{n_k-p} U^{*i}(a\dL U^{n_kl})U^{j}\del{n_k-q}=U^{*i}(a\dL U^{n_kl})U^{j}
$$
if $p=i$ and $q=j$ and it is equal $0$ otherwise. Thus
$$U^{p}\del{n_k-p} \rho_k(\sum a_{i,j}u_{n_k}^{l_{i,j}}e_{i,j}) \del{n_k-q}U^{*q}=
a_{p,q}\dL U^{n_kl_{p,q}} .$$
Hence, if $\rho_k(\sum e_{i,j}u_{n_k}^{l_{i,j}}e_{i,j})=0$, then
$a_{i,j}=0$ for all $i,j$. Thus $\rho_k$ is injective.
The fact that the map is onto follows from Lemma~\ref{Gk} . (Note
that $U^i\dL U^{*j}= U^{*(n_k-i)}\dL U^{n_k-j}=\rho_k(e_{n_k-i,n_k-j})$).

\end{proof}

\begin{theorem}\label{isom}
Let $\{n_k\}$ be an increasing sequence of positive integers with
$n_k|n_{k+1}$ (and $n_1=1$).
Then the algebra $\BDalf$ is *-isomorphic to the algebra
$C(X,A)\times_{\sigma} \mathbb{Z}$.
\end{theorem}
\begin{proof}
To construct the isomorphism we first show that, for every $k$,
\be\label{rg}
 \rho_k = \rho_{k+1} \circ \gamma_{n_k,n_{k+1}} .  \ee
Fix $k$ and write $m=n_k/n_{k+1}$.
It will suffice to apply both maps to a generating set of
$M_{n_k}(A\times_{\alpha^{n_k}}\mathbb{Z})$.
In the following computations we shall write
$\delta^{(p)}_{\mathbf{j}}$ instead of $\delta_{\mathbf{j}}$ as
before, to indicate that $\mathbf{j}=(j_1,\ldots,j_{p-1})$ is of length $p-1$
 (and $\del{j}^{(p)}$ is an element of $\cC_p$).
 For $a\in A$ we have
$$\rho_{k+1} \circ
\gamma_{n_k,n_{k+1}}(ae_{0,0})=\rho_{k+1}(\sum_{l=0}^{m-1}\alpha^{ln_k}(a)
e_{ln_k,ln_k})=$$
$$=\sum U^{*ln_k} \alpha^{n_kl}(a) \dL^{(k+1)}
U^{ln_k}=
\sum a\del{n_{k+1}-ln_k}^{(k+1)}=a\dL^{(k)}=\rho_k(ae_{0,0}). $$
Also,
$$\rho_{k+1} \circ
\gamma_{n_k,n_{k+1}}(u_{n_k}e_{0,0})=\rho_{k+1}(u_{n_{k+1}}e_{(m-1)n_k,0})+
\sum_{l=0}^{m-2} \rho_{k+1}(e_{ln_k,(l+1)n_k})=$$
$$=U^{*(m-1)n_k}\dL^{(k+1)}U^{n_{k+1}}+\sum_{l=0}^{m-2}U^{*ln_k}\dL^{(k+1)}U^{(l+1)n_k}=$$
$$=\del{n_k}^{(k+1)}U^{n_k}+\sum_{l=0}^{m-2}
\del{n_{k+1}-ln_k}^{(k+1)}U^{n_k}
=\dL^{(k)}U^{n_k}=\rho_k(u_{n_k}e_{0,0}).$$
For every $0\leq i,j \leq n_k-1$,
$$\rho_k \circ \gamma_{n_k,n_{k+1}}(e_{i,j})=\sum_{l=0}^{m-1}
\rho_k(e_{i+ln_k,j+ln_k})=\sum
U^{*(i+ln_k)}\dL^{(k+1)}U^{(j+ln_k)}=$$
$$=\sum
U^{*i}\del{n_{k+1}-ln_k}^{(k+1)}U^{j}=U^{*i}\dL^{(k)}U^{j}=\rho_k(e_{i,j}).$$
This proves (\ref{rg}) and it follows that we have a
*-homomorphism $\rho :\BDalf \\
 \rightarrow
C(X,A)\times_{\sigma}\mathbb{Z}$ whose "restriction" to
$\cB(n_k)$ is $\rho_k$. Since each map $\rho_k$ is injective, so
is $\rho$. It is left to show that $\rho$ is onto. Since the image
of $\rho_k$ is $\cG_k$, it amounts to showing that the
(increasing) union of the algebras $\cG_k$ is dense in the crossed
product. But each of the algebras $\cG_k$ contains $U$ and we know
that $C(X,A)$ is the closure of the union of the algebras $\cC_k$.
Hence the density of $\cup \cG_k$ follows.

\end{proof}
\begin{theorem}\label{simple}
Let $\{n_k\}$ be a sequence of positive integers with
$n_k|n_{k+1}$ and $n_1=1$. Then the algebra $\BDalf$ is simple if and only if,
for every $k\geq 1$, $A$ is $\alpha^{n_k}$-simple (i.e. it has no
proper closed two sided $\alpha^{n_k}$-invariant ideals).
\end{theorem}
\begin{proof}
Suppose that there is an $m\geq 1$ and an $\alpha^{n_m}$-invariant
proper ideal $J\subseteq A$. Then $J$ is also
$\alpha^{n_k}$-invariant for all $k\geq m$. For every such $k$,
the closed ideal of $A\times_{\alpha^{n_k}}\mathbb{Z}$
 generated by $J$ will be written $J_k$. Then $M_{n_k}(J_k)$ is an ideal
 (closed and proper) in $\cB(n_k)$. It is easy to check (using the
 definition of $\gamma_{n_k,n_{k+1}}$) that
$$ \gamma_{n_k,n_{k+1}}(M_{n_k}(J_k))\subseteq
M_{n_{k+1}}(J_{k+1}) $$
and, in fact,
$$ \gamma_{n_k,n_{k+1}}(M_{n_k}(J_k))=M_{n_{k+1}}(J_{k+1}) \cap
\gamma_{n_k,n_{k+1}}(M_{n_k}(A\times_{\alpha^{n_k}}\mathbb{Z})).$$
Hence $\tilde{J}:=\lim_{\rightarrow}M_{n_k}(J_k)$ is a non zero
ideal in $\BDalf$. It also follows that, for $k\geq m$,
$\tilde{J}\cap \cB(n_k)=M_{n_k}(J_k)$. Thus $\tilde{J}\neq
\BDalf$.

We now turn to prove the other direction.
We assume  that $A$ is $\alpha^{n_k}$-simple for all $k\geq 1$.
We start by showing that $C(X,A)$ is $\sigma$-simple. Let $I\subseteq C(X,A)$
 be a $\sigma$-invariant
 ideal in $C(X,A)$. Since
$$C(X,A)=\overline{ \cup \cC_k} $$
where $\cC_k$ (as defined above) form an increasing sequence of
subalgebras, it follows that
$$I=\overline{ \cup (I\cap \cC_k)} .$$
Write $I_k=I\cap \cC_k$. The algebra $\cC_k$ can be identified
with $C(\prod_{i=1}^{k-1} X_i,A)$ and $\prod_{i=1}^{k-1}X_i$ is a finite
set with $n_k$ points. Thus there are $n_k$ ideals $\{I_{k,j} :
0\leq j\leq n_k-1 \}$ in $A$ such that
$$I_k=\{ \sum_{j=0}^{n_k-1} a_j\del{j} : a_j \in I_{k,j} \} .$$
We have
$$\sigma(\sum_{j=0}^{n_k-1} a_j\del{j})=\sum
\alpha(a_j)(\del{j}\circ
\sigma_0^{-1})=\sum_{j=0}^{n_k-2}\alpha(a_j)\del{j+1}
+\alpha(a_{n_k-1})\dL .$$
Hence the $\sigma$-invariance of $I_k$ implies that
$$ \alpha(I_{k,j})\subseteq I_{k,j+1} $$
for all $0\leq j \leq n_k-2$ and
$$ \alpha(I_{k,n_k-1})\subseteq I_{k,0} .$$
Thus, for all $j$, we have
$$\alpha^{n_k}(I_{k,j})=I_{k,j} .$$
It follows from our assumptions that, for all $0\leq j \leq
n_k-1$, $I_{k,j}$ is either $A$ or $\{0\}$. In fact, the relations
above show that either $I_{k,j}=A$ for all $j$ or $I_{k,j}=\{0\}$
for all $j$. It follows that $I_k$ is either $\cC_k$ or $\{0\}$.
Since $\{I_k\}$ is an increasing sequence whose union is dense in
$I$, either $I=C(X,A)$ or $I=\{0\}$. This proves that $C(X,A)$ is
$\sigma$-simple. In order to prove that the crossed product of
$C(X,A)$ by $\sigma$ is simple it is left to show that the Connes'
spectrum of $\sigma$, $\Gamma(\sigma)$, is the full unit circle .
(See \cite[Theorem 8.11.12]{Ped}).

Let $\hat{A}$ and $\widehat{C(X,A)}$ be the set of equivalence classes
of the irreducible representations of $A$ and $C(X,A)$
respectively. For every $x\in X$ and an irreducible representation
$\tau$ of $A$ one can define an irreducible representation
$\pi=\pi_{(x,\tau)}$ of $C(X,A)$ by $\pi(f)=\tau(f(x))$ and,
conversely, every irreducible representation of $C(X,A)$ is
equivalent to some $\pi_{(x,\tau)}$. Moreover $\tau_1$ and
$\tau_2$ are equivalent if and only if the corresponding $\pi$'s
are equivalent (for the same $x$). Also, $\pi_{(x,\tau_1)}$ and
 $\pi_{(y,\tau_2)}$ are inequivalent whenever $x\neq y$.
 Hence we can write
 $$\widehat{C(X,A)}=\{\pi_{(x,\tau)} :\;x\in X,\; \tau \in
 \hat{A}\}.$$
The automorphism $\sigma$ induces a map on $\widehat{C(X,A)}$ which we
also denote by $\sigma$ (and similarly one has a map $\alpha$ on
$\hat{A}$). Then
$$\sigma(\pi_{(x,\tau)})=\pi_{(\sigma_0^{-1}(x),\tau \circ \alpha)}.$$
Suppose now that, for some $n\in \mathbb{Z}$,
$\sigma^n(\pi_{(x,\tau)})=\pi_{(x,\tau)}$. Then it follows that
$\sigma_0^n(x)=x$. But this is possible only if $n=0$ (from the definition of
$\sigma_0$) and, thus, $\sigma $ acts freely on $\widehat{C(X,A)}$. We
now use Theorem 10.4 of
\cite{OP} ( the equivalence of (i) and (v) there) to conclude that
 $\Gamma(\sigma)=\mathbb{T}$, completing
the proof.
\end{proof}

\begin{proposition}\label{inverse}
The algebras $\BDalf$ and $\cB_{\alpha^{-1}}(\{n_k\})$ are
isomorphic.
\end{proposition}
\begin{proof}
One can prove the proposition by constructing isomorphisms between
$\cB_{\alpha}(n_k)$ and $\cB_{\alpha^{-1}}(n_k)$ that are
intertwined by the connecting maps in the direct limit.
We prefer here to use the crossed product presentation. So let
$X=\prod\{0,\ldots,m_k-1\}$ (where $m_k=n_{k+1}/n_k$) and let
$\sigma_0$ be the odometer map as
above. Recall that we defined $\sigma :C(X,A)\rightarrow C(X,A)$
by $\sigma(f)(x)=\alpha(f(\sigma_0^{-1}(x)))$.
Write $\sigma'(f)(x)=\alpha^{-1}(f(\sigma_0^{-1}(x)))$. Then
$$\BDalf \cong C(X,A)\times_{\sigma}\mathbb{Z} $$ and
$$\cB_{\alpha^{-1}}(\{n_k\})\cong C(X,A)\times_{\sigma'}\mathbb{Z}
.$$
We shall write $U$ for the unitary operator that satisfies
$UfU^*=\sigma(f)$ for $f\in C(X,A)$ and such that $C(X,A)$ and $U$
generate the crossed product. Similarly we shall write $V$ for the
unitary operator that, together with $C(X,A)$, generates the other
crossed product and satisfies
$$ VfV^*=\sigma'(f) \;,\;\;\; f\in C(X,A) .$$
 Now define a map
$g$ on $X$ by
$$ g(x_1,x_2,\ldots )=(m_1-1-x_1,m_2-1-x_2,\ldots ) .$$
Then $g$ is a homeomorphism of $X$. It is also easy to check that
$$ g\circ \sigma_0=\sigma_0^{-1}\circ g .$$
Let $\Psi$ be the map from $C(X,A)\times_{\sigma}\mathbb{Z}$ into
$C(X,A)\times_{\sigma'}\mathbb{Z}$ defined by setting
$$\Psi(f)=f\circ g \;, \;\;\;f\in C(X,A) $$ and
$$\Psi(U)=V^* .$$
We have, for $f$ in $C(X,A)$,
$$V^*\Psi(f)V(x)=\sigma'^{-1}(f\circ g)(x)=\alpha(f\circ
g(\sigma_0(x)))=\alpha(f(\sigma_0^{-1}(g(x))))=$$
$$=\sigma(f)\circ
g=\Psi(\sigma(f)) .$$
Hence the map $\Psi$ is a well defined *-homomorphism on
$C(X,A)\times_{\sigma}\mathbb{Z}$. Similarly, by replacing the roles
of $\alpha$ and $\alpha^{-1}$, one can define a
*-homomorphism that is the inverse of $\Psi$, completing the proof
of the proposition.
\end{proof}

\end{section}

\begin{section}{Example : $\BDthet$}\label{th}
In this section we discuss the special case where
$A=C(\mathbb{T})$ (the continuous functions on the unit circle),
$\theta$ is a fixed irrational number
and $\alpha=\alpha_{\theta}$ is the irrational rotation by
$\theta$. Identifying $\mathbb{T}$ with $\mathbb{R}/\mathbb{Z}$ we
can write
$$ \alpha_{\theta}(f)(t)=f(t-\theta) . $$
The algebra $C(\mathbb{T})\times_{\alpha}\mathbb{Z}$ will be
written (as is costumary) $A_{\theta}$ and the resulting limit
algebra, $\BDalf$ will be denoted $\BDthet$.
It follows immediately from the results of Section~\ref{structure} that this
algebra (for a given increasing sequence $\{n_k\}$ of positive
integers each dividing the next one) is simple and has a unique
trace ( denoted $\tau$).

Given a sequence $\{n_k\}$ as above, we shall write $\delta(n_k)$
(or simply $\delta$ if the sequence is fixed)
for the associated supernatural number and $Q(\delta(n_k))$ (or
$Q(\delta)$) for the group of all rational numbers that can be
written as a quotient $m/n_k$ for some $m\in \mathbb{Z}$ and
$k\geq 1$. It is known that this group depends only on the
supernatural number of the sequence.

\begin{theorem}\label{K}
Let $\theta$ be an irrational number and $\{n_k\}$ be a sequence
of positive integers with $n_k$ dividing $n_{k+1}$ and $n_1=1$.
 Write $\delta(n_k)$
for its supernatural number and
 $\mathbb{R}^+$ for the set of all non negative real numbers.
 Then
\begin{itemize}
\item[(i)] $K_0(\BDthet)$ is isomorphic, as an ordered group, to
$(\qnk+\theta \mathbb{Z}, \\ (\qnk+\theta \mathbb{Z})\cap
\mathbb{R}^+)$.
\item[(ii)] Let $\tau_*$ be the map on $K_0(\BDthet)$ induced by
the unique trace. Then
$$\tau_*(K_0(\BDthet))=\qnk+\theta \mathbb{Z}.$$
\end{itemize}
\end{theorem}
\begin{proof}
Since $\cB(n_k)$ is the algebra of $n_k\times n_k$ matrices over
the irrational rotation algebra $A_{n_k\theta}$, the ordered group
$K_0(\cB(n_k))$ is order isomorphic to $(\mathbb{Z}+\theta
\mathbb{Z},(\mathbb{Z}+n_k\theta \mathbb{Z})\cap \mathbb{R}^+)$.
Let $\tau_0$ be the tracial state on $C(\mathbb{T})$ that
one gets by integrating with respect to the normalized Lebegue
measure (and then the unique trace $\tau$ on $\BDthet$ is
$\Phi_1(\tau_0)$ in the notation of Proposition~\ref{trace}) then
we write (as in Proposition~\ref{trace}) $\tau_k=\tau_0 \circ
\cE_k$ (where $\cE_k$ is the conditional expectation from
$A_{n_k\theta}$ to $C(\mathbb{T})$) and $\tilde{\tau}_k$ for the
induced trace on $\cB(n_k)$.
Let $j:A_{n_k\theta}\rightarrow \cB(n_k)=M_{n_k}(A_{n_k\theta})$
be the map $j(X)=Xe_{0,0}$. Then $j_*$ is an isomorphism of the
$K_0$ groups and, since $\tilde{\tau}_k\circ j=(1/n_k)\tau$, we
have
$$ (\tilde{\tau}_k)_*\circ j_*=\frac{1}{n_k}(\tau_k)_* .$$
Recall that we denote by $u_{n_k}$ the unitary that, together with
a copy of $C(\mathbb{T})$ generates the crossed product
$C(\mathbb{T})\times_{\alpha_{\theta n_k}}\mathbb{Z}=A_{\theta
n_k}$ . Thus $u_{n_k}$ can be identified with the operator of
rotation by $\theta n_k$ on $L^2(\mathbb{T})$.

 It is known that $(\tau_k)_*$ is an
order isomorphism from $K_0(A_{\theta n_k})$ onto
$\mathbb{Z}+\theta n_k \mathbb{Z}$ (\cite[10.11.6]{Bl} or \cite[Example 5.8]
{RLL}). It follows that $(\tilde{\tau}_k)_*$ is an isomorphism of
$K_0(\cB(n_k))$ onto $(1/n_k)\mathbb{Z}+\theta \mathbb{Z}$.
Write $\gamma_k$ for the map $\gamma_{n_k,n_{k+1}}$ and $(\gamma_k)_*$
for the map it induces on the $K_0$ groups. Hence
$$(\gamma_k)_*: K_0(\cB(n_k)) \rightarrow K_0(\cB(n_{k+1})) .$$

We have $\tilde{\tau}_{k+1}\circ \gamma_k =\tilde{\tau}_k$
(Proposition~\ref{trace}) and, thus, $(\tilde{\tau}_{k+1})_*\circ
(\gamma_k)_*=(\tilde{\tau}_k)_* $ and the map
$$ \lim_{\rightarrow}(\tilde{\tau}_k)_* :
\lim_{\rightarrow}K_0(\cB(n_k)) \rightarrow \mathbb{R} $$
is a well defined order isomorphism into $\mathbb{R}$. In fact,
since $\tau=\lim_{\rightarrow}\tilde{\tau}_k$
(Proposition~\ref{trace}), the map
$\lim_{\rightarrow}(\tilde{\tau}_k)_*$ is $\tau_*$. We also know
that $\lim_{\rightarrow}K_0(\cB(n_k))$ is isomorphic to
$K_0(\BDthet)$ as ordered groups (\cite[Theorem 6.3.2]{RLL}).
 It follows that
$$\tau_*:K_0(\BDthet) \rightarrow \mathbb{R} $$
is an order isomorphism. Its image is
$$\cup ((1/n_k)\mathbb{Z}+\theta \mathbb{Z}) =
Q(\delta(n_k))+\theta \mathbb{Z} .$$

\end{proof}

\vspace{3mm}

\begin{corollary}\label{isomorphism}
We have $\cB_{\theta_1}(\{n_k\}) \cong \cB_{\theta_2}(\{m_k\}) $
if and only if $\delta(n_k)=\delta(m_k)$ and either
$\theta_1+\theta_2$ or $\theta_1-\theta_2$ lies in
$Q(\delta(n_k))$ ($=Q(\delta(m_k))$).
\end{corollary}
\begin{proof}
Assume first that the condition on the supernatural numbers and the
 $\theta$'s holds. It follows from Proposition~\ref{inverse} that
 we can assume that $\theta_1-\theta_2$ lies in $Q(\delta (n_k))$
 (replacing $\theta_2$ by $-\theta_2$ if necessary).
Then, for some $k\geq 1$, $n_k\theta_1-n_k\theta_2$ lies in
$\mathbb{Z}$. In fact, $n_l\theta_1-n_l\theta_2 \in \mathbb{Z}$
for every $l\geq k$. The algebras $A_{n_l\theta_1}$ and
$A_{n_l\theta_2}$ are then isomorphic. In fact, if we write
$\psi_l$ for this isomorphism and apply $\psi_l$ entrywise we
get an isomorphism $\tilde{\psi}_l$ of $\cB_{\theta_1}(n_l)$ onto
$\cB_{\theta_2}(n_l)$. These isomorphisms intertwines the
connecting maps and we get an isomorphism of the limit algebras.

For the other direction, assume the two algebras are isomorphic.
Write $\eta$ for the isomorphism $\eta : \cB_{\theta_1}(\{n_k\})
\rightarrow \cB_{\theta_2}(\{m_k\})$. If $\tau$ is the unique
tracial state on $ \cB_{\theta_2}(\{m_k\})$ then $\tau \circ \eta
$ is the unique tracial state on $ \cB_{\theta_1}(\{n_k\})$. Using
Proposition~\ref{K} we have
$$Q(\delta(n_k))+\theta_1 \mathbb{Z} = (\tau \circ \eta)_*(K_0(
 \cB_{\theta_1}(\{n_k\})))=\tau_*(\eta_*(K_0(
 \cB_{\theta_1}(\{n_k\}))))=$$
 $$=\tau_*(K_0(\cB_{\theta_2}(\{m_k\})))=
 Q(\delta(m_k))+\theta_2 \mathbb{Z} .$$
It follows that, for every $k \geq 1$, there are $l \geq 1$ and
$p,c$ in $\mathbb{Z}$ such that $1/n_k=p/m_l +c\theta_2$. Since
$\theta_2$ is irrational, $c=0$ and $n_k|m_l$. This shows that
$\delta(n_k)$ divides $\delta(m_k)$ and , by symmetry, they are
equal. Hence
$$Q(\delta(n_k))+\theta_1 \mathbb{Z} =Q(\delta(n_k))+\theta_2 \mathbb{Z} .$$
Thus there is some $k\geq 1$ and integers $a,b,c,d$ such that
$\theta_1=b/n_k+a\theta_2 $ and $\theta_2=c/n_k+d\theta_1$.
Combining these equalities we find that $ad=1$ . Hence either
$a=d=1$ (and then $\theta_1-\theta_2 \in Q(\delta(n_k))$) or $a=d=-1$ (and
then $\theta_1+\theta_2 \in Q(\delta(n_k))$).

\end{proof}

\vspace{3mm}

One can also compute the $K_1$ group of the algebra $\BDthet$
using the continuity of $K_1$.

\begin{proposition}\label{K1}
For every $\theta$ and $\{n_k\}$ as above,
$$K_1(\BDthet) \cong Q(\delta(n_k)) \oplus \mathbb{Z} .$$
\end{proposition}
\begin{proof}
It is known that $K_1(A_{n_k\theta}) \cong \mathbb{Z}\oplus\mathbb{Z} $
for every $k\geq 1$ and the generators of the group are the class
$[u_{n_k}]$ of the unitary $u_{n_k}$ and the class $[v]$ where $v$
 is the function in $C(\mathbb{T})$ defined by $v(z)=z$ (\cite[Example VIII.5.2]
 {D} and \cite[10.11.6]{Bl}) .
 Fix $k\geq 1$, write $m=n_{k+1}/n_k$ and $\gamma_k$ for the map
$\gamma_{n_k,n_{k+1}}$. The $K_1$ group of
$M_{n_k}(A_{n_k\theta})$ is again $\mathbb{Z}\oplus\mathbb{Z}$ and the
generators are the classes of $\tilde{u}_k :=
u_{n_k}e_{0,0}+\sum_{j=1}^{n_k-1}e_{j,j} $ and $\tilde{v}_k
:=ve_{0,0}+\sum_{j=1}^{n_k-1}e_{j,j} $.
We compute
$$\gamma_k(\tilde{u}_k)=u_{n_{k+1}}e_{n_k(m-1),0}+\sum_{l=0}^{m-2}
e_{ln_k,(l+1)n_k}+\sum_i e_{i,i} $$
and
$$\gamma_k(\tilde{v}_k)=\sum_{l=0}^{m-1}\alpha_{ln_k\theta}(v)e_{ln_k,ln_k}
+\sum_i e_{i,i} $$
where the last sums in both equations run over all integers $0\leq
i \leq n_{k+1}-1$ that are not multiples of $n_k$.
But then it follows easily that
$$(\gamma_k)_*([\tilde{u}_k])=[\tilde{u}_{k+1}] $$
and
$$(\gamma_k)_*([\tilde{v}_k])=m[\tilde{v}_{k+1}] .$$
Hence, viewing $(\gamma_k)_*$ as a map from
$\mathbb{Z}\oplus \mathbb{Z}$ to $\mathbb{Z}\oplus \mathbb{Z}$, we have
$$(\gamma_k)_* :(a,b) \mapsto (n_{k+1}/n_k a,b) .$$
Using the continuity of $K_1$ (\cite[Proposition 8.2.7]{RLL}) we
get the result.

\end{proof}

\vspace{3mm}

We close with the observation that every matrix algebra over an
algebra of this class is again an algebra of this class. The
precise statement is presented in the following proposition.

\begin{proposition}
For every $p \in \mathbb{N}$,
$$M_p(\BDthet) \cong \cB_{\theta/p}(\{pn_k\}) .$$
Hence, for $p >1$, the algebras $\BDthet$ and $M_p(\BDthet)$ are
non isomorphic.
\end{proposition}
\begin{proof}
Write
$$\BDthet=\lim_{\rightarrow}(\cB_{\theta}(n_k),\gamma_k) $$
where  $\gamma_k$ stands for $\gamma_{n_k,n_{k+1}}$.
Then
$$M_p(\BDthet)=\lim_{\rightarrow}(M_p(\cB_{\theta}(n_k)),(\gamma_k)_p)=
\lim_{\rightarrow}(M_p(M_{n_k}(A_{\theta n_k})),(\gamma_k)_p) $$
where $(\gamma_k)_p$ is defined by applying $\gamma_k$ entrywise.
We also have,
$$\cB_{\theta/p}(\{pn_k\})=\lim_{\rightarrow}(\cB_{\theta/p}(pn_k),\gamma'_k)
=\lim_{\rightarrow}(M_{pn_k}(A_{(\theta/p)pn_k}),\gamma'_k) $$
where $\gamma'_k$ is $\gamma_{pn_k,pn_{k+1}}$.
An element of $M_p(M_{n_k}(A_{\theta n_k}))$ can be written, in an
obvious way, as an $pn_k\times pn_k$ matrix over $A_{\theta n_k}$.
Now perform the canonical shuffle on the matrix by applying the
permutation
$$t_k: (0, \ldots ,pn_k-1) \mapsto
(0,p,2p,\ldots,(n_k-1)p,1,p+1,\ldots, p(n_k-1)-1,pn_k-1) $$
on the rows and columns of the matrix. This defines an isomorphism
$$\psi_k : M_p(M_{n_k}(A_{\theta n_k})) \rightarrow
M_{pn_k}(A_{\theta n_k}) $$
that satisfies
$$ \psi_{k+1}\circ (\gamma_k)_p=\gamma'_k \circ \psi_k \;\;\;,k\geq
1.$$
The map $\psi:=\lim \psi_k$ is the required isomorphism.

The last statement follows from Corollary~\ref{isomorphism} since
$\theta - \theta /p $ is irrational.
\end{proof}
\end{section}

\end{document}